\documentclass{amsart}
\usepackage[T1]{fontenc}
\usepackage{color}
\usepackage{textcomp}
\usepackage{subfig}
\usepackage{amsmath}
\usepackage{amssymb}
\usepackage{amsthm}
\usepackage[english]{babel}
\usepackage{graphicx}
\usepackage{tikz}
\usepackage{hyperref}
\usepackage{enumitem}

\theoremstyle{definition}
\newtheorem{definition}{Definition}
\newtheorem{theorem}{Theorem}
\newtheorem{lemm}[theorem]{Lemma}

\newtheorem{conjecture}[theorem]{Conjecture}

\begin{document}
\title[Non-Uniquely Ergodic Interval Exchange Transformations with Flips]{Existence of a Non-Uniquely Ergodic Interval Exchange Transformation with Flips Possessing Three Invariant Measures}
\author {Aleksei Kobzev}
\address{Faculty of Mathematics, National Research University Higher School of Economics, Usacheva St. 6, 119048 Moscow, Russia}
\email{akobzev@hse.ru}

\begin{abstract}
We present the first explicit example of an interval exchange transformation with flips (FIET) possessing three distinct invariant ergodic measures. The proof is based on a generalization of M. Keane's method, using the Rauzy induction adapted for FIETs, which contributes to the study of the ergodic properties of this class of dynamical systems.
\end{abstract}
\maketitle
\section{Introduction}
\subsection {Overview of the main results}
The main object of our research is interval exchange transformations with flips (FIETs). An FIET is a piecewise isometric map of an interval to itself with a finite number of jump discontinuities, where the orientation of at least one of the intervals of continuity is reversed, this property will be referred to as a flip (for more details on FIETs, see \cite{Ng87}). FIETs are a generalization of classic interval exchange transformations (IETs).

IETs arise as the first return map to a transversal for an orientable measured foliation on an orientable surface or, for instance, as a first return map to a transversal for a billiard flow in a rational polygon.
IETs were first mentioned by V. Arnold in 1963 in his work \cite{Arn}, where he considered IETs on three intervals.

IETs are very actively studied. The general case was first described by V. Oseledets in 1968 in \cite{Odc}. M. Keane showed in \cite{Kn75} that a typical IET is minimal, meaning the orbit of every point is dense. Later, H. Masur and W. Veech showed in \cite{M82} and \cite{V82} that almost all IETs are uniquely ergodic. Also, in a work by A. Katok \cite{Kt}, an upper bound on the number of invariant ergodic measures was obtained: an IET cannot have more than $m$ invariant ergodic measures, where $m$ is the genus of a surface with the property that the given IET is the first return map of the vertical flow to a horizontal transversal, i.e., $\lfloor \frac{n}{2} \rfloor$, where $n$ is the number of intervals. The first proof of the effectiveness of this bound was obtained by E. Sataev \cite{St}. This result was also obtained by J. Fickenscher \cite{Fk}.

The dynamics of FIETs differs from the dynamics of classical IETs. For instance, in \cite{Ng89}, A. Nogueira proved that a typical FIET is not minimal. Due to this, the class of minimal FIETs remains largely unexplored. In particular, much less is known about the ergodic properties of minimal FIETs. Nevertheless, some progress has been made. The Hausdorff dimension of the set of minimal FIETs was estimated by A. Skripchenko and S. Troubetzkoy in \cite{ST18}. Also, in \cite{BL23}, A. Linero Bas and G. Soler López established the existence of a minimal non-uniquely ergodic example of a FIET on 10 intervals with two distinct invariant ergodic measures. Subsequently, in \cite{BL25}, they provided an example on 6 intervals with two distinct invariant ergodic measures.

It is worth noting that A. Skripchenko and S. Troubetzkoy \cite{ST18} described a general procedure to construct minimal non-uniquely ergodic FIETs with an arbitrary number of invariant measures. Their method consists of extending a minimal IET with $n-1$ subintervals on the interval $I$ by adding a small additional interval. Under the constructed map, one subinterval of $I$ is mapped onto the additional interval with an orientation reversal, while the additional interval is mapped back into $I$, also with a reversal, such that the first return map to $I$ coincides with the original IET. Although this construction yields minimal examples, they arise essentially as extensions of standard IETs without flips and are non-generic. The novelty of our work lies in constructing an explicit irreducible example that does not rely on such a method. In this paper, we present a proof of the following result:

\begin{theorem}\label{zero} 
There exists an irreducible minimal non-uniquely ergodic FIET with three distinct invariant ergodic measures.
\end{theorem}

The proof is based on a method proposed by M. Keane in \cite{Kn77} for classical IETs, which involves estimating the measures of certain subintervals for each of the invariant measures and proving their difference. We generalize this approach to the case of FIETs. The key tool of analysis is a generalized Rauzy induction, adapted for FIETs, the structure of which inherits the main ideas of the original construction by G. Rauzy \cite{Rz}.

\subsection{Applications of FIETs} 
The reason for the interest in studying FIETs is their connection to other areas of mathematics, particularly to translation surfaces and billiards in rational polygons. C. Danthony and A. Nogueira in \cite{DN90} exploited linear involutions with flips, which are a generalization of FIETs, to show that almost all measured foliations on non-orientable surfaces have a compact leaf. In \cite{HR18}, P. Hubert and O. Paris-Romaskevich showed a connection between tiling billiards and FIETs, which allowed them to prove several theorems concerning triangle tiling billiards.

\subsection{Organization of the paper}
The paper is organized as follows. In Section \ref{sec:fiet}, we recall the basic definitions of classical IETs and formally introduce FIETs. We also describe the generalized Rauzy induction algorithm adapted for FIETs, which serves as our main analytical tool. 

Section \ref{sec:main} is devoted to the proof of the main result. In Section \ref{subsec:construction}, we define the specific combinatorial data of our example and construct a closed Rauzy path, calculating the corresponding transition matrix. Section \ref{subsec:estimates} contains a sequence of four technical lemmas establishing precise upper and lower bounds for the lengths of specific subintervals with respect to three limits of the transition matrices. Specifically, the first two lemmas estimate the distributions of lengths for two prospective measures, while the last two lemmas secure the bounds for the third one. In Section \ref{subsec:non-unique}, we combine these estimates to prove Theorem \ref{final} (which immediately implies Theorem \ref{zero}), showing that the three measures are pairwise distinct due to the incompatibility of their mass distributions. Finally, Section \ref{subsec:conclusion} provides an upper bound on the Hausdorff dimension of the parameter set corresponding to our example and formulates an open conjecture.

\newpage

\subsection{Acknowledgments}
I am deeply grateful to A. Skripchenko for suggesting the problem and for her constant guidance and encouragement throughout the preparation of this paper. I would also like to thank S. Troubetzkoy for providing important comments. The research leading to these results has received funding from the Basic Research Program at HSE University (HSE-BR-2025-62).

\subsection{Some notations}

\begin{description}[itemsep=1pt, topsep=3pt, labelwidth=0.6cm, leftmargin=!]
    \item[$\Theta$] The transition matrix corresponding to the Rauzy graph.

    \item[$x^{(m)}$] A vector $x^{(m)} = (x^{(m)}_1, \dots, x^{(m)}_n)$ at step $m$, with components satisfying $\sum_i x^{(m)}_i = 1$.

    \item[$|\mathbf{v}|$] The $L_1$-norm of a vector $\mathbf{v}$, i.e., $|\mathbf{v}| = \sum_i |v_i|$.

    \item[$\overline{\Theta}$] The normalized operator, $\overline{\Theta} = \frac{\Theta}{|\Theta|}$. When applied to a vector, it implies $\overline{\Theta}x = \frac{\Theta x}{|\Theta x|}$.

    \item[$\mu_q(\cdot)$] Notation for measures. It will be shown in the paper that each such measure is an invariant ergodic probability measure and that measures with different indices are distinct. An important notational convention is:
    $$
    x^{(m)} = \mu_q(I^{(m)}), \quad x^{(m)}_i = \frac{\lambda_q(I^{(m)}_i)}{\lambda_q(I^{(m)})}.
    $$
    Here it is understood that before passing to $x^{(m)}$, the specific measure $\mu_q$ is fixed and agreed upon. When only one measure is involved, we typically use the index $q$. When two distinct measures are compared, we use indices $i$, $j$.

    \item[$I^{(m)}, I^{(m)}_i$] Subintervals. Specifically, $I^{(m)}_i$ denotes the $i$-th subinterval at the $m$-th step of the first return map.
\end{description}




\section{Interval Exchange Transformations with Flips}\label{sec:fiet}

In this section, we will present the formal definitions and tools that we will use to prove the results presented in the current paper.

\subsection{Basic definitions}\label{app:basic}

\begin{definition}
An IET is a map $T: I \to I$ determined by a pair $(\lambda, \pi)$, where:
\begin{itemize}
    \item $\lambda = (\lambda_1, \dots, \lambda_n) \in \mathbb{R}_+^n$ is a vector of lengths. Let $L = \sum_{i=1}^n \lambda_i$.
    \item $\pi = (\pi_0, \pi_1) \in \mathcal{S}_n \times \mathcal{S}_n$ is a pair of permutations on $\{1, \dots, n\}$.
\end{itemize}
The permutations define partitions of the domain and range of $T$. The interval $I$ is partitioned into $n$ consecutive subintervals, whose orderings are given by $\pi_0$ for the domain and $\pi_1$ for the range. For each label $k \in \{1, \dots, n\}$, there is a unique subinterval of length $\lambda_k$ in both the domain and the range partition.

The map $T$ sends the domain subinterval with label $k$ to the range subinterval with the same label $k$ by a translation. At the discontinuity points, $T$ is defined to be right‑continuous.
\end{definition}

$$
\begin{tikzpicture}[scale=0.8]
\draw[ultra thick, red] (0,0) -- (1,0);
\draw[ultra thick, green] (1,0) -- (6,0);
\draw[ultra thick, blue] (6,0) -- (10,0);
\node at (0.5,0.5) {\color{black} $\lambda_1$};
\node at (3.5,0.5) {\color{black} $\lambda_2$};
\node at (8,0.5) {\color{black} $\lambda_3$};
\draw[ultra thick, red] (9,-0.5) -- (10,-0.5);
\draw[ultra thick, green] (0,-0.5) -- (5,-0.5);
\draw[ultra thick, blue] (5,-0.5) -- (9,-0.5);
\node at (9.5,-1) {\color{black} $\lambda_1$};
\node at (2.5,-1) {\color{black} $\lambda_2$};
\node at (7,-1) {\color{black} $\lambda_3$};
\end{tikzpicture}
$$
\begin{center}
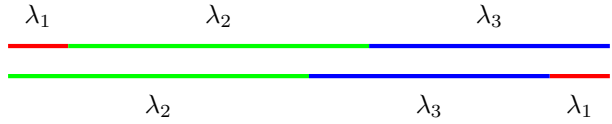
\textbf{Figure 1:} An IET for $n=3$ with permutations $\pi_0 = (1, 2, 3)$ and $\pi_1 = (2, 3, 1)$.
\end{center}

\begin{definition}
  An IET $T$ is called an FIET if the partitioned interval not only has its subintervals rearranged but also has their orientation reversed. The map $T$ is defined by a triple $(\lambda, \pi, F)$:
\begin{itemize}
    \item $(\lambda, \pi)$ are as in the definition of an IET.
    \item $F \subseteq \{1, \dots, n\}$ is a set of indices called the \emph{flip set}. It is the set of intervals that reverse their orientation under the map $T$.
\end{itemize}
The discontinuity points of $T$ are the endpoints of the subintervals $I_i$, this set is finite and therefore of Lebesgue measure zero. At each discontinuity point, $T$ is defined to be right‑continuous.
\end{definition}

$$
\begin{tikzpicture}[scale=0.8]
\draw[ultra thick, red, ->] (0,0) -- (1,0);
\draw[ultra thick, green, ->] (1,0) -- (6,0);
\draw[ultra thick, blue, ->] (6,0) -- (10,0);
\node at (0.5,0.5) {\color{black} $\lambda_1$};
\node at (3.5,0.5) {\color{black} $\lambda_2$};
\node at (8,0.5) {\color{black} $\lambda_3$};
\draw[ultra thick, red, <-] (9,-0.5) -- (10,-0.5);
\draw[ultra thick, green, ->] (0,-0.5) -- (5,-0.5);
\draw[ultra thick, blue, <-] (5,-0.5) -- (9,-0.5);
\node at (9.5,-1) {\color{black} $\lambda_1$};
\node at (2.5,-1) {\color{black} $\lambda_2$};
\node at (7,-1) {\color{black} $\lambda_3$};
\end{tikzpicture}
$$
\begin{center}
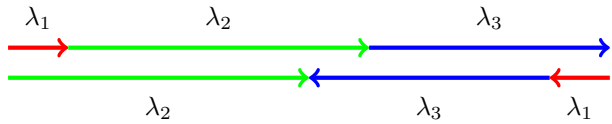
\textbf{Figure 2:} An FIET for $n=3$ with permutations $\pi_0 = (1, 2, 3)$, $\pi_1 = (2, 3, 1)$ and $F = \{1, 3\}$.
\end{center}

\subsection{Rauzy Induction} The main tool used in the proof of the paper's results is a generalization of Rauzy induction.
\begin{definition}
Rauzy induction is the main renormalization procedure that allows us to construct a new IET from an initial one, which has the exact same orbit structure but is defined on a smaller support interval. Rauzy induction was first introduced by Gerard Rauzy for IETs in his 1979 paper \cite{Rz}.
\newline
 Right Rauzy induction:
Consider an IET $T: I \to I$ on $n$ subintervals of $I = [0, 1)$. We take the rightmost subintervals in the domain and range of the map. Then we induce $T$ on $[0, 1 - \lambda)$, where $\lambda = \min(\lambda_{a_n}, \lambda_{a_i})$ is the length of the shorter of the two rightmost subintervals of $I$ and $TI$. We will say that the longer interval is the "winner" ($W$) and the shorter one is the "loser" ($L$).

We adapt Rauzy induction for FIETs, it will act slightly differently on our defined pair of permutations:
The induction action transforms $(\pi_0, \pi_1, F)$ into $(\pi'_0, \pi'_1, F')$. \\

\paragraph{Case (a): $\boldsymbol{\lambda_{\pi_0(n)} > \lambda_{\pi_1(n)}}$:} 
$ $ \\
In this case, $W = \pi_0(n)$, $L = \pi_1(n)$, and $\pi'_0 = \pi_0$. Only $\pi_1$ changes.
We find the position of the Winner in the permutation being modified: $k = \pi_1^{-1}(W) = \pi_1^{-1}(\pi_0(n))$.
\begin{enumerate}
\item[a.1)] $W \notin F$, insert $L$ immediately after $W$, $\pi'_0 = \pi_0$, $F' = F \Delta \{\pi_1(n)\}$.
$$
\pi'_1(i) =
\begin{cases}
\pi_1(i) & \text{if } 1 \le i \le k \\
\pi_1(n) & \text{if } i = k+1 \\
\pi_1(i-1) & \text{if } k+2 \le i \le n
\end{cases}
$$

\item[a.2)] $W \in F$ 
, insert $L$ immediately before $W$, $\pi'_0 = \pi_0$, $F' = F \Delta \{\pi_1(n)\}$.
$$
\pi'_1(i) =
\begin{cases}
\pi_1(i) & \text{if } 1 \le i \le k-1 \\
\pi_1(n) & \text{if } i = k \\
\pi_1(i-1) & \text{if } k+1 \le i \le n
\end{cases}
$$
\end{enumerate}

\paragraph{Case (b): $\boldsymbol{\lambda_{\pi_1(n)} > \lambda_{\pi_0(n)}}$:} 
$ $ \\
In this case, $W = \pi_1(n)$, $L = \pi_0(n)$, and $\pi'_1 = \pi_1$. Only $\pi_0$ changes.
We find the position of the Winner in the permutation being modified: $k = \pi_0^{-1}(W) = \pi_0^{-1}(\pi_1(n))$.
\begin{enumerate}
\item[b.1)] $W \notin F$, insert $L$ immediately after $W$, $\pi'_1 = \pi_1$, $F' = F \Delta \{\pi_0(n)\}$.
$$
\pi'_0(i) =
\begin{cases}
\pi_0(i) & \text{if } 1 \le i \le k \\
\pi_0(n) & \text{if } i = k+1 \\
\pi_0(i-1) & \text{if } k+2 \le i \le n
\end{cases}
$$

\item[b.2)] $W \in F$, insert $L$ immediately before $W$, $\pi'_1 = \pi_1$, $F' = F \Delta \{\pi_0(n)\}$.
$$
\pi'_0(i) =
\begin{cases}
\pi_0(i) & \text{if } 1 \le i \le k-1 \\
\pi_0(n) & \text{if } i = k \\
\pi_0(i-1) & \text{if } k+1 \le i \le n
\end{cases}
$$
\end{enumerate}
\end{definition} 

These four systems completely and formally describe our new, more complex induction operation. We will denote the induction as $R_{\gamma}$, where $\gamma$ is a sequence of induction steps.

Let's write down the 4 possible induction steps in a diagram:

b.1) If $a_n \notin F$, then the set $F$ does not change, and the diagram is as follows:

\begin{tikzpicture}[
    ->,
    >=stealth,
    shorten >=5pt,
    auto,
    node distance=7cm,
    semithick
]

\tikzset{
    state/.style={fill=white, draw=none, text=black, rounded corners, inner sep=4pt}
}

\node[state] (A) {$\underbrace{\begin{Bmatrix}
a_1&...&&&...&a_n\\
...&...&a_n&a_j&...&a_i
\end{Bmatrix}}_{\pi}$};
\node[state] (B) [right of=A] {$\underbrace{\begin{Bmatrix}
a_1&...&&&...&a_n\\
...&a_n&a_i&a_j&...&...
\end{Bmatrix}}_{\pi'}$};

\path[->] (A) edge node {b} (B);

\end{tikzpicture}

b.2) If $a_n \in F$, then the set $F$ becomes $F \Delta \{a_i\}$, and the diagram is as follows:

\begin{tikzpicture}[
    ->,
    >=stealth,
    shorten >=5pt,
    auto,
    node distance=7cm,
    semithick
]

\tikzset{
    state/.style={fill=white, draw=none, text=black, rounded corners, inner sep=4pt}
}

\node[state] (A) {$\underbrace{\begin{Bmatrix}
a_1&...&&&...&a_n\\
...&...&a_n&a_j&...&a_i
\end{Bmatrix}}_{\pi}$};
\node[state] (B) [right of=A] {$\underbrace{\begin{Bmatrix}
a_1&...&&&...&a_n\\
...&a_i&a_n&a_j&...&...
\end{Bmatrix}}_{\pi'}$};

\path[->] (A) edge node {b} (B);

\end{tikzpicture}

a.1) If $a_i \notin F$, then the set $F$ does not change, and the diagram is as follows:

\begin{tikzpicture}[
    ->,
    >=stealth,
    shorten >=5pt,
    auto,
    node distance=7cm,
    semithick
]

\tikzset{
    state/.style={fill=white, draw=none, text=black, rounded corners, inner sep=4pt}
}

\node[state] (A) {$\underbrace{\begin{Bmatrix}
a_1&...&a_i&a_{i+1}&...&a_n\\
...&...&&&...&a_i
\end{Bmatrix}}_{\pi}$};
\node[state] (B) [right of=A] {$\underbrace{\begin{Bmatrix}
a_1&...&a_i&a_n&a_{i+1}&...\\
...&...&&&...&a_i
\end{Bmatrix}}_{\pi'}$};

\path[->] (A) edge node {a} (B);

\end{tikzpicture}

a.2) If $a_i \in F$, then the set $F$ becomes $F \Delta \{a_n\}$, and the diagram is as follows:

\begin{tikzpicture}[
    ->,
    >=stealth,
    shorten >=5pt,
    auto,
    node distance=7cm,
    semithick
]

\tikzset{
    state/.style={fill=white, draw=none, text=black, rounded corners, inner sep=4pt}
}

\node[state] (A) {$\underbrace{\begin{Bmatrix}
a_1&...&a_i&a_{i+1}&...&a_n\\
...&...&&&...&a_i
\end{Bmatrix}}_{\pi}$};
\node[state] (B) [right of=A] {$\underbrace{\begin{Bmatrix}
a_1&...&a_n&a_i&a_{i+1}&...\\
...&...&&&...&a_i
\end{Bmatrix}}_{\pi'}$};

\path[->] (A) edge node {a} (B);

\end{tikzpicture}

Also, figures illustrating the generalized Rauzy induction are presented below, divided into two groups: with and without changes.

$$
\begin{tikzpicture}[scale=0.8]

\draw[ultra thick, red, ->] (0,0) -- (1,0);
\draw[ultra thick, blue, ->] (1,0) -- (3,0);
\draw[ultra thick, green, ->] (3,0) -- (7,0);

\draw[ultra thick, red, ->] (9,0) -- (13,0);
\draw[ultra thick, blue, ->] (13,0) -- (15,0);
\draw[ultra thick, green, ->] (15,0) -- (16,0);

\draw[ultra thick, green, ->] (0,-0.6) -- (4,-0.6);
\draw[ultra thick, blue, ->] (4,-0.6) -- (6,-0.6);
\draw[ultra thick, red, ->] (6,-0.6) -- (7,-0.6);

\draw[ultra thick, green, ->] (9,-0.6) -- (10,-0.6);
\draw[ultra thick, blue, ->] (10,-0.6) -- (12,-0.6);
\draw[ultra thick, red, ->] (12,-0.6) -- (16,-0.6);

\draw[ultra thick, red, ->] (0,0-1.7) -- (1,0-1.7);
\draw[ultra thick, blue, ->] (1,0-1.7) -- (3,0-1.7);
\draw[ultra thick, green, ->] (3,0-1.7) -- (6,0-1.7);

\draw[ultra thick, red, ->] (9,0-1.7) -- (12,0-1.7);
\draw[ultra thick, blue, ->] (13,0-1.7) -- (15,0-1.7);
\draw[ultra thick, green, ->] (12,0-1.7) -- (13,0-1.7);

\draw[ultra thick, green, ->] (0,-0.6-1.7) -- (3,-0.6-1.7);
\draw[ultra thick, blue, ->] (4,-0.6-1.7) -- (6,-0.6-1.7);
\draw[ultra thick, red, ->] (3,-0.6-1.7) -- (4,-0.6-1.7);

\draw[ultra thick, green, ->] (9,-0.6-1.7) -- (10,-0.6-1.7);
\draw[ultra thick, blue, ->] (10,-0.6-1.7) -- (12,-0.6-1.7);
\draw[ultra thick, red, ->] (12,-0.6-1.7) -- (15,-0.6-1.7);

\draw[dashed] (6,-3) -- (6,0.3);
\draw[dashed] (15,-3) -- (15,0.3);

\draw[ultra thick, ->, black] (3,-0.8) -- (3,-1.3) node[midway, right] {b.1};
\draw[ultra thick, ->, black] (12,-0.8) -- (12,-1.3) node[midway, right] {a.1};
\end{tikzpicture}
$$
\begin{center}
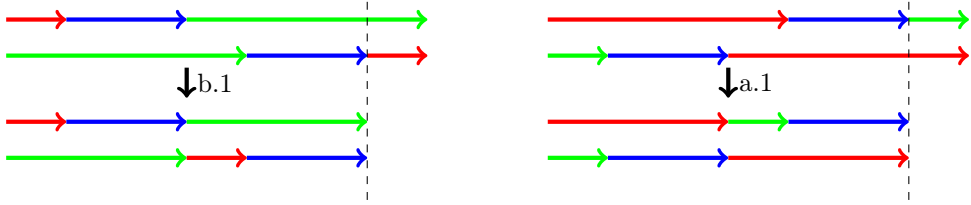
\textbf{Figure 3:} Steps when F does not change.
\end{center}

$$
\begin{tikzpicture}[scale=0.8]

\draw[ultra thick, red, ->] (0,0) -- (1,0);
\draw[ultra thick, blue, ->] (1,0) -- (3,0);
\draw[ultra thick, green, ->] (3,0) -- (7,0);

\draw[ultra thick, red, ->] (9,0) -- (13,0);
\draw[ultra thick, blue, ->] (13,0) -- (15,0);
\draw[ultra thick, green, ->] (15,0) -- (16,0);

\draw[ultra thick, green, <-] (0,-0.6) -- (4,-0.6);
\draw[ultra thick, blue, ->] (4,-0.6) -- (6,-0.6);
\draw[ultra thick, red, ->] (6,-0.6) -- (7,-0.6);

\draw[ultra thick, green, ->] (9,-0.6) -- (10,-0.6);
\draw[ultra thick, blue, ->] (10,-0.6) -- (12,-0.6);
\draw[ultra thick, red, <-] (12,-0.6) -- (16,-0.6);

\draw[ultra thick, red, ->] (0,0-1.7) -- (1,0-1.7);
\draw[ultra thick, blue, ->] (1,0-1.7) -- (3,0-1.7);
\draw[ultra thick, green, ->] (3,0-1.7) -- (6,0-1.7);

\draw[ultra thick, red, ->] (10,0-1.7) -- (13,0-1.7);
\draw[ultra thick, blue, ->] (13,0-1.7) -- (15,0-1.7);
\draw[ultra thick, green, ->] (9,0-1.7) -- (10,0-1.7);

\draw[ultra thick, green, <-] (1,-0.6-1.7) -- (4,-0.6-1.7);
\draw[ultra thick, blue, ->] (4,-0.6-1.7) -- (6,-0.6-1.7);
\draw[ultra thick, red, <-] (0,-0.6-1.7) -- (1,-0.6-1.7);

\draw[ultra thick, green, <-] (9,-0.6-1.7) -- (10,-0.6-1.7);
\draw[ultra thick, blue, ->] (10,-0.6-1.7) -- (12,-0.6-1.7);
\draw[ultra thick, red, <-] (12,-0.6-1.7) -- (15,-0.6-1.7);

\draw[dashed] (6,-3) -- (6,0.3);
\draw[dashed] (15,-3) -- (15,0.3);

\draw[ultra thick, ->, black] (3,-0.8) -- (3,-1.3) node[midway, right] {b.2};
\draw[ultra thick, ->, black] (12,-0.8) -- (12,-1.3) node[midway, right] {a.2};
\end{tikzpicture}
$$
\begin{center}
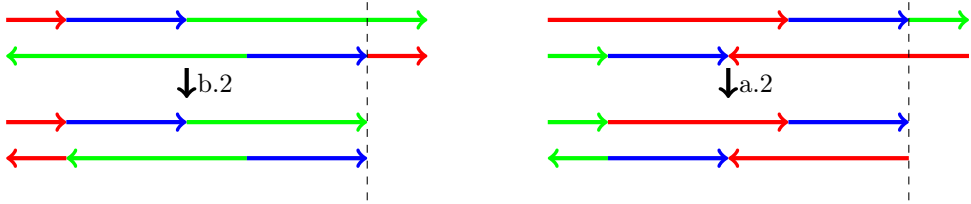
\textbf{Figure 4:} Steps when F changes.
\end{center}
$ $

Now we will define the transition matrix for each of these actions. Note that the old lengths are expressed in terms of the new ones according to the following principle: all interval lengths remain the same, except for the original winner's interval, whose length becomes the sum of the new winner's and loser's interval lengths. The matrix for each action will look as follows:

$$
E_{i,j} = 
\begin{cases} 
1, & \text{if } i = j \ \text{or}\ (i = w \ \text{and}\ j = l,\ l \neq w), \\
0, & \text{otherwise}
\end{cases}
$$

\section{Main Result}\label{sec:main}

\subsection{Construction of the non-uniquely ergodic example}\label{subsec:construction}
Consider the following permutation and set $F$, which define a FIET:
$$ \pi = 
\begin{Bmatrix}
\pi_0\\
\pi_1
\end{Bmatrix} =
\begin{Bmatrix}
1&2&3&4&5&6&7&8\\
4&5&6&7&2&1&8&3
\end{Bmatrix}, \ 
F = \{2, 3, 4, 5, 6, 7\}
$$
\\
Let $p_1, p_2, p_3 > 0$ and define the following Rauzy paths:

$$\gamma_1 = aaaabbaaaaabbbaaab,$$
$$\gamma_2(p_1) = a^{p_1},$$
$$\gamma_3 = baaabaaaaabbb,$$
$$\gamma_4(p_2) = a^{p_2},$$
$$\gamma_5 = baababa,$$
$$\gamma_6(p_3) = b^{p_3}.$$
Then the total path is:
$$
\gamma_p = \gamma_1\gamma_2(p_1)\gamma_3\gamma_4(p_2)\gamma_5\gamma_6(p_3)
$$
$$
\gamma_p = aaaabbaaaaabbbaaaba^{p_1}baaabaaaaabbba^{p_2}baababab^{p_3}
$$
Following this path on the Rauzy diagram, we arrive at the following pair $(\pi, F)$:
$$
R_{\gamma_p}((\pi, F)) = (\pi', F'),
$$
where
$$
\pi' = 
\begin{Bmatrix}
\pi_0\\
\pi_1
\end{Bmatrix} =
\begin{Bmatrix}
1&4&2&3&5&6&7&8\\
3&5&6&7&4&1&8&2
\end{Bmatrix}, \ 
F' = \{2, 3, 4, 5, 6, 7\}
$$

Note that we can consider $R^{ind_{2,3,4}}_{\gamma_p}$ - the induced action of the Rauzy induction, where we observe how the induction acts on intervals 2, 3, and 4.

Note that
$$
R^{ind_{2,3,4}}_{\gamma_p}(
\begin{Bmatrix}
2&3&4\\
4&2&3
\end{Bmatrix}
) =
\begin{Bmatrix}
4&2&3\\
3&4&2
\end{Bmatrix}
$$
This is a cyclic permutation with period 3:
$$
\begin{tikzpicture}[
    ->,
    >=stealth,
    shorten >=5pt,
    auto,
    node distance=5cm,
    semithick
]

\tikzset{
    state/.style={fill=white, draw=none, text=black, rounded corners, inner sep=4pt}
}

\node[state] (A) {$\begin{Bmatrix}
2&3&4\\
4&2&3
\end{Bmatrix}$};
\node[state] (B) [right of=A] {$\begin{Bmatrix}
4&2&3\\
3&4&2
\end{Bmatrix}$};
\node[state] (C) [right of=B] {$\begin{Bmatrix}
3&4&2\\
2&3&4
\end{Bmatrix}$};

\path[->] (A) edge node {$R^{ind_{2,3,4}}_{\gamma_p}$} (B);
\path[->] (B) edge node {$R^{ind_{2,3,4}}_{\gamma_p}$} (C);
\path[->] (C) edge [bend right=50] node {$R^{ind_{2,3,4}}_{\gamma_p}$} (A);

\end{tikzpicture}
$$

Thus, we find that the full cycle for our initial FIET $(\pi, F)$ is:
$$
\tilde{\gamma_p} = \gamma_p \gamma_p \gamma_p,
$$

$$
R_{\tilde{\gamma_p}}((\pi, F)) = (\pi, F).
$$

Let's compute the transition matrix. We are interested in the transition matrix for $\gamma_p$, as it defines the dynamics of the exchange, and $\tilde{\gamma_p}$ is a power of this matrix.

Now let's compute the transition matrices for each component of the path $\gamma_p$:

\begin{align*}
\gamma_1 = \begin{pmatrix}
1 & 1 & 1 & 2 & 1 & 2 & 1 & 1 \\
0 & 1 & 0 & 0 & 0 & 0 & 0 & 0 \\
1 & 2 & 1 & 1 & 1 & 2 & 2 & 2 \\
1 & 1 & 1 & 1 & 0 & 1 & 1 & 1 \\
0 & 0 & 0 & 0 & 1 & 0 & 0 & 0 \\
0 & 1 & 0 & 1 & 0 & 1 & 0 & 1 \\
0 & 0 & 0 & 0 & 0 & 0 & 1 & 0 \\
1 & 1 & 0 & 0 & 0 & 0 & 1 & 1
\end{pmatrix},
&&
\gamma_2(p_1) = \begin{pmatrix}
1 & 0 & 0 & 0 & 0 & 0 & 0 & 0 \\
0 & 1 & 0 & 0 & 0 & 0 & 0 & 0 \\
0 & 0 & 1 & 0 & 0 & 0 & 0 & 0 \\
0 & 0 & 0 & 1 & 0 & 0 & 0 & 0 \\
0 & 0 & p_1 & 0 & 1 & 0 & 0 & 0 \\
0 & 0 & 0 & 0 & 0 & 1 & 0 & 0 \\
0 & 0 & 0 & 0 & 0 & 0 & 1 & 0 \\
0 & 0 & 0 & 0 & 0 & 0 & 0 & 1
\end{pmatrix}.
\end{align*}

\begin{align*}
\gamma_3 = \begin{pmatrix}
1 & 1 & 1 & 1 & 1 & 0 & 0 & 1 \\
0 & 1 & 0 & 1 & 0 & 0 & 0 & 0 \\
0 & 0 & 1 & 0 & 1 & 0 & 0 & 0 \\
0 & 0 & 1 & 1 & 1 & 0 & 0 & 1 \\
0 & 0 & 0 & 0 & 1 & 0 & 0 & 0 \\
1 & 1 & 0 & 1 & 0 & 1 & 0 & 0 \\
0 & 0 & 0 & 0 & 0 & 0 & 1 & 0 \\
0 & 0 & 0 & 0 & 0 & 0 & 0 & 1
\end{pmatrix},
&&
\gamma_4(p_2) = \begin{pmatrix}
1 & 0 & 0 & 0 & 0 & 0 & 0 & 0 \\
0 & 1 & 0 & 0 & 0 & 0 & 0 & 0 \\
0 & 0 & 1 & 0 & 0 & 0 & 0 & 0 \\
0 & 0 & 0 & 1 & 0 & 0 & 0 & 0 \\
0 & 0 & 0 & 0 & 1 & 0 & 0 & 0 \\
0 & 0 & 0 & 0 & 0 & 1 & 0 & 0 \\
0 & 0 & 0 & 0 & 0 & p_2 & 1 & 0 \\
0 & 0 & 0 & 0 & 0 & 0 & 0 & 1
\end{pmatrix}.
\end{align*}

\begin{align*}
\gamma_5 = \begin{pmatrix}
1 & 1 & 0 & 1 & 0 & 0 & 0 & 0 \\
0 & 1 & 0 & 1 & 0 & 0 & 0 & 1 \\
0 & 0 & 1 & 0 & 0 & 0 & 0 & 0 \\
1 & 0 & 0 & 1 & 0 & 0 & 0 & 0 \\
0 & 0 & 0 & 0 & 1 & 0 & 0 & 0 \\
0 & 0 & 0 & 0 & 0 & 1 & 1 & 0 \\
0 & 0 & 0 & 0 & 0 & 0 & 1 & 0 \\
0 & 0 & 0 & 0 & 0 & 1 & 1 & 1
\end{pmatrix},
&&
\gamma_6(p_3) = \begin{pmatrix}
1 & 0 & 0 & 0 & 0 & 0 & 0 & 0 \\
0 & 1 & 0 & 0 & 0 & 0 & 0 & 0 \\
0 & 0 & 1 & 0 & 0 & 0 & 0 & 0 \\
0 & 0 & 0 & 1 & 0 & 0 & 0 & 0 \\
0 & 0 & 0 & 0 & 1 & 0 & 0 & 0 \\
0 & 0 & 0 & 0 & 0 & 1 & 0 & 0 \\
0 & 0 & 0 & 0 & 0 & 0 & 1 & 0 \\
0 & p_3 & 0 & 0 & 0 & 0 & 0 & 1
\end{pmatrix}.
\end{align*}

Then the resulting matrix is\footnote{Some computations are available at \href{https://github.com/Kobzetsu/FIET_with_three_measures/blob/main/FIET_with_three_measures_calc.ipynb}{link}}:

\begin{align*}
\scalebox{0.9}{
$\Theta_{\gamma_{p}} = $
$\begin{pmatrix}
9 & 8p_3 + 7 & p_1 + 4 & 13 & p_1 + 5 & p_2 + 6 & p_2 + 7 & 8 \\
1 & p_3 + 1 & 0 & 2 & 0 & 0 & 0 & 1 \\
9 & 9p_3 + 8 & p_1 + 3 & 14 & p_1 + 4 & 2p_2 + 6 & 2p_2 + 8 & 9 \\
6 & 6p_3 + 5 & 3 & 9 & 3 & p_2 + 4 & p_2 + 5 & 6 \\
0 & 0 & p_1 & 0 & p_1 + 1 & 0 & 0 & 0 \\
4 & 4p_3 + 3 & 1 & 6 & 1 & 3 & 3 & 4 \\
0 & 0 & 0 & 0 & 0 & p_2 & p_2 + 1 & 0 \\
3 & 4p_3 + 3 & 1 & 5 & 1 & p_2 + 2 & p_2 + 3 & 4
\end{pmatrix}$
}
\end{align*}

In the next section, we will introduce a sequence of parameters $\{p^i\}_{i=1}^{\infty} = \{p_1^{(i)}, p_2^{(i)}, p_3^{(i)}\}_{i=1}^{\infty}$ and, based on constraints on this sequence, prove some estimates for the subintervals of each of the measures $\mu_2, \mu_5, \mu_7$, where $\mu_q$ corresponds to the vector of partition interval lengths as 
$$\lim_{m \to +\infty} \overline{\Theta}_{1}\overline{\Theta}_{2}...\overline{\Theta}_{m} e_q, \text{ where } \Theta_i = \Theta_{\gamma_{p^{3i-2}}}\Theta_{\gamma_{p^{3i-1}}}\Theta_{\gamma_{p^{3i}}}.$$

\subsection{Proof of the estimates}\label{subsec:estimates}
\begin{lemm}
The following estimates for measure $\mu_7$ hold for every step $m$:
\begin{enumerate}
\item $x^{(m)}_7 > x^{(m)}_5$;
\item $2x^{(m)}_7 > x^{(m)}_8$;
\item $2x^{(m)}_7 > x^{(m)}_4$;
\item $2x^{(m)}_7 > x^{(m)}_1$;
\item $4x^{(m)}_7 > x^{(m)}_3$;
\item $x^{(m)}_7 > \frac{1}{7}$;
\item $x^{(m)}_2 < \frac{1}{p^{(m)}_1}$.
\item $x^{(m)}_5 < \frac{1}{10}$.
\end{enumerate}
\end{lemm}

\begin{proof}
The base case is verified by constructing inequalities over the first base vector corresponding to $e_7$. We trivially calculate $$\Theta_{\gamma_{p}} e_7= (p_2+7, 0, 2p_2+8, p_2+5, 0, 3, p_2 +1, p_2+3)^T,$$ structurally confirming all inequalities (for $p_2 > 20$). \\
From the construction it is clear that $x^{(m)}_6 > x^{(m)}_2$ and $x^{(m)}_3 > x^{(m)}_7$.
Assume that for step $m+1$ the required estimates hold.
By the induction hypothesis we have
\begin{align*}
|\Theta_m x^{(m+1)}| = &\ 32x^{(m+1)}_1 + (32p^{(m)}_3+27)x^{(m+1)}_2 + (3p^{(m)}_1+12)x^{(m+1)}_3 + 49x^{(m+1)}_4 \\
&+ (3p^{(m)}_1+15)x^{(m+1)}_5 + (6p^{(m)}_2+21)x^{(m+1)}_6 + (6p^{(m)}_2+27)x^{(m+1)}_7 + 32x^{(m+1)}_8 \\
&> \frac{1}{2}\,p^{(m)}_2.
\end{align*}

\noindent\textbf{Proof of $x^{(m)}_7 > x^{(m)}_5$.} Consider the difference
\begin{align*}
|\Theta_m x^{(m+1)}|\,(x^{(m)}_7 - x^{(m)}_5)
&= p^{(m)}_2 x^{(m+1)}_6 + (p^{(m)}_2+1)x^{(m+1)}_7 - p^{(m)}_1 x^{(m+1)}_3 - (p^{(m)}_1+1)x^{(m+1)}_5 \\
&> p^{(m)}_2 x^{(m+1)}_6 + p^{(m)}_2 x^{(m+1)}_7 - p^{(m)}_1 x^{(m+1)}_3 - p^{(m)}_1 x^{(m+1)}_5 \\
&> p^{(m)}_2 x^{(m+1)}_6 + (d-1)p^{(m)}_1 x^{(m+1)}_7 - p^{(m)}_1 x^{(m+1)}_3.
\end{align*}
Using $4x^{(m+1)}_7 > x^{(m+1)}_3$ and $d-1 > 4$, we obtain $$(d-1)p^{(m)}_1 x^{(m+1)}_7 - p^{(m)}_1 x^{(m+1)}_3 > 0.$$

\noindent\textbf{Proof of $2x^{(m)}_7 > x^{(m)}_8$.} Consider
\begin{align*}
|\Theta_m x^{(m+1)}|\,(2x^{(m)}_7 - x^{(m)}_8)
&= (p^{(m)}_2-2)x^{(m+1)}_6 + (p^{(m)}_2-1)x^{(m+1)}_7 - 3x^{(m+1)}_1 \\
&\quad  - (4p^{(m)}_3+3)x^{(m+1)}_2 - x^{(m+1)}_3 - 5x^{(m+1)}_4 - x^{(m+1)}_5 - 4x^{(m+1)}_8 \\
&> (p^{(m)}_2-2)x^{(m+1)}_6 + (p^{(m)}_2-1)x^{(m+1)}_7 - \frac{5}{d} - 3x^{(m+1)}_1 \\
&\quad - x^{(m+1)}_3 - 5x^{(m+1)}_4 - x^{(m+1)}_5 - 4x^{(m+1)}_8.
\end{align*}
For $p^{(m)}_2$ sufficiently large ($p^{(m)}_2-1 > 30$) this expression is positive.
The estimates $2x^{(m)}_7 > x^{(m)}_4$ (requires $p^{(m)}_2-3 > 58$, $d>56$) and $2x^{(m)}_7 > x^{(m)}_1$ (requires $d>72$) are proved similarly.

\noindent\textbf{Proof of $4x^{(m)}_7 > x^{(m)}_3$.} This follows directly from the inequalities: 
$$2x^{(m+1)}_7 > x^{(m+1)}_1 \text{ and } 2x^{(m+1)}_7 > x^{(m+1)}_8.$$

\noindent\textbf{Proof of $x^{(m)}_7 > \frac{1}{7}$.} We need to show
\[
7(\Theta_m x^{(m+1)})_7 - |\Theta_m x^{(m+1)}| > 0,
\]
i.e.
\begin{multline*}
(p^{(m)}_2-21)x^{(m+1)}_6 + (p^{(m)}_2-20)x^{(m+1)}_7 - \frac{33}{d} - 32x^{(m+1)}_1 \\
- (3p^{(m)}_1+12)x^{(m+1)}_3 - 49x^{(m+1)}_4 - (3p^{(m)}_1+15)x^{(m+1)}_5 - 32x^{(m+1)}_8 > 0,
\end{multline*}
which holds for $d > 33\cdot 7$.

\noindent\textbf{Proof of $x^{(m)}_2 < \frac{1}{p^{(m)}_1}$.} Since $|\Theta_m x^{(m+1)}| > 2p^{(m)}_1$ and $x^{(m+1)}_3 > x^{(m+1)}_4$, we obtain
\[
x^{(m)}_2 < \frac{x^{(m+1)}_1 + x^{(m+1)}_3 + x^{(m+1)}_4 + x^{(m+1)}_8 + \frac{p^{(m)}_3+1}{p^{(m+1)}_1}}{2p^{(m)}_1}
< \frac{1+\frac{2}{d}}{2p^{(m)}_1} < \frac{1}{p^{(m)}_1}.
\]

\noindent\textbf{Proof of $x^{(m)}_5 < \frac{1}{10}$.} Using $|\Theta_m x^{(m+1)}| > \frac{1}{2}p^{(m)}_2$ and the previously established inequalities,
\[
\begin{aligned}
x^{(m)}_5 &< \frac{p^{(m)}_1 x^{(m+1)}_3 + (p^{(m)}_1+1)x^{(m+1)}_5}{\frac{1}{2}p^{(m)}_2} < \frac{4p^{(m)}_1 x^{(m+1)}_7 + (p^{(m)}_1+1)x^{(m+1)}_7}{\frac{1}{2}p^{(m)}_2} \\
&< \frac{5p^{(m)}_1+1}{\frac{d}{2}p^{(m)}_1} < \frac{11}{d}.
\end{aligned}
\]
Thus for $d > 110$ we have $x^{(m)}_5 < \frac{1}{10}$.
\end{proof}
\begin{lemm}
The following estimates for $\mu_5$ hold for every step $m$:
\begin{enumerate}
\item $x^{(m)}_2 < \frac{1}{p^{(m)}_1}$;
\item $x^{(m)}_6 < \frac{7}{p^{(m)}_1}$;
\item $x^{(m)}_7 < \frac{1}{p^{(m)}_1}$;
\item $x^{(m)}_5 > \frac{1}{4}$;
\item $2x^{(m)}_3 > x^{(m)}_1$;
\item $x^{(m)}_3 + x^{(m)}_5 > x^{(m)}_1$;
\item $3x^{(m)}_6 + x^{(m)}_7 > x^{(m)}_4$;
\item $x^{(m)}_6 + x^{(m)}_7 \ge x^{(m)}_8$.
\end{enumerate}
\end{lemm}

\begin{proof}
The base case is verified by constructing inequalities over the first base vector corresponding to $e_5$. We trivially calculate $$\Theta_{\gamma_{p}} e_5= (p_1+5, 0, p_1+4, 3, p_1+1, 1, 0, 1)^T,$$ structurally confirming all inequalities (for $p_1 > 11$). \\

Note that inequalities (5)--(8) follow immediately from a direct component-wise comparison of the corresponding rows of the transition matrix $\Theta_{\gamma_p}$, and thus hold unconditionally for any non-negative vector $x^{(m+1)}$.\\

Assume that for step $m+1$ the required estimates hold.
By the induction hypothesis we have
$$|\Theta_m x^{(m+1)}| > p^{(m)}_1.$$

\noindent\textbf{Proof of $x^{(m)}_2 < \frac{1}{p^{(m)}_1}$.}  
This is proved similarly to the case of $\mu_7$.

\noindent\textbf{Proof of $x^{(m)}_6 < \frac{7}{p^{(m)}_1}$.}  
\[
\begin{aligned}
x^{(m)}_6 &< \frac{4x^{(m+1)}_1+\frac{5}{d}+x^{(m+1)}_3+6x^{(m+1)}_4+x^{(m+1)}_5+3x^{(m+1)}_6+3x^{(m+1)}_7+4x^{(m+1)}_8}{p^{(m)}_1} \\
&< \frac{6 + \frac{5}{d}}{p^{(m)}_1} < \frac{7}{p^{(m)}_1}.
\end{aligned}
\]

\noindent\textbf{Proof of $x^{(m)}_7 < \frac{1}{p^{(m)}_1}$.}  
\[
x^{(m)}_7 < \frac{p^{(m)}_2 x^{(m+1)}_6 + (p^{(m)}_2+1)x^{(m+1)}_7}{p^{(m)}_1}
< \frac{\frac{7}{d^2}+\frac{2}{d^2}}{p^{(m)}_1} < \frac{1}{p^{(m)}_1}.
\]



\noindent\textbf{Proof of $x^{(m)}_5 > \frac{1}{4}$.}  
It suffices to show
\[
4p^{(m)}_1 x^{(m+1)}_3 + 4(p^{(m)}_1+1)x^{(m+1)}_5 > |\Theta_m x^{(m+1)}|,
\]
which is equivalent to
\begin{multline*}
(p^{(m)}_1-12)x^{(m+1)}_3 + (p^{(m)}_1-11)x^{(m+1)}_5 - 32x^{(m+1)}_1 - (32p^{(m)}_3+27)x^{(m+1)}_2 \\
- 49x^{(m+1)}_4 - (6p^{(m)}_2+21)x^{(m+1)}_6 - (6p^{(m)}_2+27)x^{(m+1)}_7 - 32x^{(m+1)}_8 > 0.
\end{multline*}
This inequality holds for $p^{(m)}_1 > 45$.
\end{proof}
\begin{lemm}
If $c>10$, then for measure $\mu_2$ for every step $m$ and for each index $i \in \{1,3,4,5,6,7,8\}$ the inequality
\[
c\,x^{(m)}_2 > x^{(m)}_i
\]
holds.
\end{lemm}

\begin{proof}
The base case is obvious. Assume that for step $m+1$ the required inequality holds.
The proof for $m$ is based on comparing the coefficients of the components of the vector $x^{(m+1)}$ in the expressions for $x^{(m)}_k$. We compare the coefficient of $x^{(m+1)}_2$ in $x^{(m)}_2$, multiplied by $c$, with the sum of coefficients in the other rows. The coefficient for $x^{(m)}_2$ is $p^{(m)}_3+1$. Below are the sums of coefficients for the other coordinates:

\[
\begin{aligned}
x^{(m)}_1 &: \ 8p^{(m)}_3+2p^{(m)}_1+2p^{(m)}_2+59,\\
x^{(m)}_3 &: \ 9p^{(m)}_3+2p^{(m)}_1+4p^{(m)}_2+61,\\
x^{(m)}_4 &: \ 6p^{(m)}_3+2p^{(m)}_2+41,\\
x^{(m)}_5 &: \ 2p^{(m)}_1+1,\\
x^{(m)}_6 &: \ 4p^{(m)}_3+25,\\
x^{(m)}_7 &: \ 2p^{(m)}_2+1,\\
x^{(m)}_8 &: \ 4p^{(m)}_3+2p^{(m)}_2+22.
\end{aligned}
\]

The largest coefficient is in the row for $x^{(m)}_3$. Since $p^{(m)}_3 > 2p^{(m)}_1+4p^{(m)}_2+61$, it is sufficient to take $c>10$ for the inequality
\[
c(p^{(m)}_3+1) > 9p^{(m)}_3+2p^{(m)}_1+4p^{(m)}_2+61
\]
to hold. Consequently $c\,x^{(m)}_2 > x^{(m)}_i$ for every $i \neq 2$.
\end{proof}

\begin{lemm}
If $b>33$, then for $\mu_2$ for every step $m$ the inequality
\[
x^{(m)}_2 > \frac{1}{b}
\]
holds.
\end{lemm}

\begin{proof}
The base case is obvious. Assume that for step $m+1$ the required inequality holds.
We need to prove that $b(\Theta_m x^{(m+1)})_2 - |\Theta_m x^{(m+1)}| > 0$. After substitution and simplification we obtain
\[
\begin{aligned}
& (b-32)p^{(m)}_3 x^{(m+1)}_2 - 3p^{(m)}_1\bigl(x^{(m+1)}_3 + x^{(m+1)}_5\bigr) - 6p^{(m)}_2\bigl(x^{(m+1)}_6 + x^{(m+1)}_7\bigr) - 49 \\
& \quad > (b-32)p^{(m)}_3 x^{(m+1)}_2 - c\bigl(6p^{(m)}_1 + 12p^{(m)}_2\bigr)x^{(m+1)}_2 - 49 \\
& \quad = \bigl((b-32)p^{(m)}_3 - 6c(p^{(m)}_1 + 2p^{(m)}_2)\bigr)x^{(m+1)}_2 - 49 \\
& \quad > (b-33)p^{(m)}_3 x^{(m+1)}_2 - 49 > (b-33)p^{(m)}_3 \cdot \frac{1}{b} - 49.
\end{aligned}
\]
The last inequality is equivalent to $(b-33)(p^{(m)}_3-49) > 33\cdot 49$.\\ Since $p^{(m)}_3 > d^3 > 50^3$, it suffices to take $b>33$.
\end{proof}
\subsection{Non-unique ergodicity}\label{subsec:non-unique}
\begin{theorem}\label{final}
    Let there be a sequence of parameters $\{p_1^{(i)}, p_2^{(i)}, p_3^{(i)}\}_{i=1}^{\infty}$, where $p_1^{(i+1)} = dp_3^{(i)} = d^2p_2^{(i)} = d^3p_1^{(i)}$ and $p_1^{(1)} > d > 50^3$. \\
    The matrix $\Theta_i = \Theta_{\gamma_{p^{3i-2}}}\Theta_{\gamma_{p^{3i-1}}}\Theta_{\gamma_{p^{3i}}}$ is the transition matrix for $\tilde{\gamma_p}$.
    Let a FIET be given as $(\pi, \lambda, F)$, where $\lambda \in \bigcap_{m=1}^{\infty} \overline{\Theta}_{1}...\overline{\Theta}_{m} v$, $v \in \Delta_7$, and $\pi$ and $F$ are from the beginning of the section \ref{app:basic}.
    Then $\mu_2$, $\mu_5$, and $\mu_7$ are pairwise distinct invariant ergodic measures, where $\mu_q$ corresponds to the vector of partition interval lengths as 
$$\lim_{m \to +\infty} \overline{\Theta}_{1}\overline{\Theta}_{2}...\overline{\Theta}_{m} e_q.$$
\end{theorem}

\begin{proof}
Consider the initial partition $\{I_i\}_{i=1}^8$ of the interval $I$ with length vector $\lambda = (\lambda_1,\dots,\lambda_8)$. Under the given constraints on the parameters $\{p_1^{(i)}, p_2^{(i)}, p_3^{(i)}\}_{i=1}^{\infty}$, the lemmas proved above imply the following estimates on the frequencies of visits to the subintervals for each measure:

\begin{itemize}
\item For $\mu_7$: $\mu_7(I_7) > 1/7$, while $\mu_7(I_2) < 1/p_1$ and $\mu_7(I_5) < 1/10$.
\item For $\mu_5$: $\mu_5(I_5) > 1/4$, while $\mu_5(I_2) < 1/p_1$ and $\mu_5(I_7) < 1/p_1$.
\item For $\mu_2$: $\mu_2(I_2) > 1/34$.
\end{itemize}

Now compare $\mu_5$ and $\mu_7$. If they were equal, they would have to assign the same measure to every measurable set. However,
\[
\mu_5(I \setminus I_7) + \mu_7(I_7) > \left(1 - \frac{1}{p_1}\right) + \frac{1}{7} > 1
\]
for sufficiently large $p_1$. This is impossible because the left‑hand side is at most $\mu_5(I \setminus I_7) + \mu_5(I_7) = 1$ if the measures coincide. Hence $\mu_5 \neq \mu_7$.

Similarly,
\[
\mu_7(I \setminus I_5) + \mu_5(I_5) > \left(1 - \frac{1}{10}\right) + \frac{1}{4} > 1,
\]
which again contradicts equality.

For $\mu_2$ and $\mu_7$ we have
\[
\mu_7(I \setminus I_2) + \mu_2(I_2) > \left(1 - \frac{1}{p_1}\right) + \frac{1}{34} > 1,
\]
so $\mu_2 \neq \mu_7$. The difference between $\mu_2$ and $\mu_5$ is proved analogously. Thus the three measures are pairwise distinct.
\end{proof}

\begin{proof}[Proof of Theorem \ref{zero}]
Theorem \ref{zero} follows directly from Theorem \ref{final}.
\end{proof}

\subsection{Conclusion}\label{subsec:conclusion}

A recent work by C. Fougeron~\cite{CF} introduced a methodology for obtaining upper bounds on the Hausdorff dimension of parameter sets arising from Rauzy–Veech induction. 
Let $\Delta_{\mathrm{NUE}}$ denote the set of parameters (in the sense of Theorem~\ref{final}) that remain forever in the strongly connected component of the Rauzy diagram that contains our constructed example. By Proposition~3.29 of~\cite{CF}, this component is \emph{non‑extremal} (i.e., it has at least one outgoing edge).  
Equivalently,
\[
\Delta_{\mathrm{NUE}} = \bigcap_{m=1}^{\infty} \overline{\Theta}_{1}\cdots\overline{\Theta}_{m} v,
\qquad v \in \Delta_7.
\]
Applying Theorem~B' of~\cite{CF} to this non‑extremal component, we obtain an upper bound for the Hausdorff dimension. For a component with $|\mathcal{A}|$ labels, the theorem gives
\[
\dim_H \Delta(F) \le |\mathcal{A}|-3 + \frac{h}{|\mathcal{A}|-1} < |\mathcal{A}|-2,
\]
where $0<h<|\mathcal{A}|-1$ is the entropy of the suspension flow. In our case $|\mathcal{A}|=8$, and consequently
\[
\dim_H(\Delta_{\mathrm{NUE}}) < 6.
\]
This shows that the set of parameters giving rise to our example has Hausdorff dimension strictly smaller than the dimension of the full parameter space.
This inequality shows that the family of examples we have constructed is not typical in the parameter space: its Hausdorff dimension is strictly smaller than the dimension of the full space. This constitutes a significant step toward a general quantitative understanding of how the Hausdorff dimension of the parameter set depends on the number of distinct ergodic invariant measures possessed by a FIET.

In connection with the question of the possible number of ergodic invariant measures for FIETs, we formulate the following conjecture.

\begin{conjecture}
A FIET on $n$ intervals can have at most $\left\lfloor \frac{n}{2} \right\rfloor - 1$ distinct invariant ergodic measures.
\end{conjecture}

\end{document}